\documentclass[11pt, reqno]{amsart}
\usepackage{color}
\usepackage[utf8]{inputenc}
\usepackage[T1]{fontenc}

\usepackage{amsmath}
\usepackage{amssymb}
\usepackage{amsfonts}
\usepackage{mathrsfs}
\usepackage{amsbsy}
\usepackage{relsize}
\usepackage[colorlinks, linkcolor=red, citecolor=blue, urlcolor=blue, pagebackref, hypertexnames=false]{hyperref}
\usepackage[hyperpageref]{backref}
\usepackage{enumitem}

\usepackage{tikz}

\numberwithin{equation}{section}
\newcounter{bbb}
\numberwithin{bbb}{section}

\newtheorem{theorem}[bbb]{Theorem}
\newtheorem{lemma}[bbb]{Lemma}
\newtheorem{proposition}[bbb]{Proposition}

\newtheorem{definition}[bbb]{Definition}
\newtheorem{remark}[bbb]{Remark}

\renewcommand{\geq}{\geqslant}

\renewcommand{\leq}{\leqslant}

\usepackage{color}

\newcommand{\RN}{\mathbb{R}^N}
\newcommand{\N}{\mathbb{N}}
\newcommand{\R}{\ensuremath{\mathbb{R}}}

\title[Hardy-H\'{e}non equation]{On the Fujita exponent for a Hardy-H\'{e}non equation with a spatial-temporal forcing term}
\author[M. Majdoub]{Mohamed Majdoub}
\address[M. Majdoub]{Department of Mathematics, College of Science, Imam Abdulrahman Bin Faisal University, P. O. Box 1982, Dammam, Saudi Arabia.\newline Basic and Applied Scientific Research Center, Imam Abdulrahman Bin Faisal University, P.O. Box 1982, 31441, Dammam, Saudi Arabia}
\email{\sl mmajdoub@iau.edu.sa}
\email{\sl med.majdoub@gmail.com}

\begin{document}

\begin{abstract}
The purpose of this work is to analyze the wellposedness and the blow-up of solutions of the higher-order parabolic semilinear equation
\[
    u_t+(-\Delta)^{d}u=|x|^{\alpha}|u|^{p}+\zeta(t){\mathbf w}(x) \ \quad\mbox{for }
    (x,t)\in\mathbb{R}^{N}\times(0,\infty),
\]
where $d\in (0,1)\cup \mathbb{N}$, $p>1$, $-\alpha\in(0,\min(2d,N))$ or $\alpha\geq 0$ and $\zeta$ as well as ${\mathbf w}$ are suitable given functions. Given $p\geq \frac{N-2d\sigma+\alpha}{N-2d\sigma-2d}$ and setting $p_c=\frac{N(p-1)}{2d+\alpha}$, $\ell=\frac{N p_c}{N+2(\sigma+1)d p_c}$, we prove that for any data $u_0\in L^{p_c,\infty}(\mathbb{R}^N)$ and $\textbf{w}\in L^{\ell,\infty}(\mathbb{R}^N)$ with small norms there exists a unique global-in-time solution under the hypotheses  $\zeta(t)=t^{\sigma}$, $\sigma\in (-1,0)$ and $N>2d$ in the space $C_{b}([0,\infty);L^{p_c,\infty}(\mathbb{R}^N))$. As a by-product, small Lebesgue data global existence follows and in particular, unconditional uniqueness holds in $C_{b}([0,\infty);L^{p_c}(\mathbb{R}^N))$ provided $p\in (\frac{N+\alpha}{N-2d},\infty)$. If either $m\in (-\infty,0]$ and $p\in (1,\frac{N-2dm+\alpha}{N-2dm-2d})$ or $m>0$ and $p>1$ where $\zeta(t)=O(t^m)$, $t\rightarrow\infty$ ($m\in \mathbb{R}$), then all solutions blow up under the additional condition $\int_{\mathbb{R}^N}\textbf{w}(x)\,dx>0$. As a consequence, we deduce that the corresponding Fujita critical exponent is a function of $\sigma$ and reads $p_{F}(\sigma)=\frac{N-2d\sigma+\alpha}{N-2d\sigma-2d}$ if $-1<\sigma<0$ and infinity otherwise.
\end{abstract}
\date{\today}
\maketitle

$$\rule[0.05cm]{10cm}{0.05cm}$$


\section{Introduction and main results}

In this paper we are interested in the study of the critical behaviour of global in time solutions of the initial value problem
\begin{equation}\label{eq:base}
\begin{cases}
u_t+(-\Delta)^{d}u=|x|^{\alpha}|u|^{p}+\zeta(t){\mathbf w}(x) \ \quad\mbox{for }
    (x,t)\in\mathbb{R}^{N}\times(0,\infty),\\
u(x,0)=u_0(x),\quad x\in \mathbb{R}^{N},
\end{cases}
\end{equation}
where $N\geq 1$, $\alpha\in\R$, $p>1$ and $\zeta, {\mathbf w}$ are given functions. The operator $(-\Delta)^{d}$ for $d\geq 1$ integer is the $d$-fold composition of the Laplacian with itself while for $d\in (0,1)$, $(-\Delta)^{d}$ should be interpreted as the fractional Laplacian defined by
$(-\Delta)^d\,u:=\mathcal{F}^{-1}\left(|\xi|^{2d}\mathcal{F}(u)\right)$ where $\mathcal{F}$ and $\mathcal{F}^{-1}$ are the Fourier transform and  inverse Fourier transform, respectively \cite{M}. We will consider the situation $\alpha\geq 0$ in which case $x\in \mathbb{R}^{N}$ and $0<-\alpha<\min(2d,N)$ with $x\in \mathbb{R}^{N}\setminus \{0\}$.

Reaction-diffusion equations with forcing term arises in many physical phenomena and biological species theories, such as the concentration of diffusion of some fluid, the density of some biological species, and heat conduction phenomena, see \cite{Hu, GV, QS, MT1, MT2} and references therein.

Our main goal in this paper is to find the critical exponent which separates the existence and nonexistence of global solutions of \eqref{eq:base}.
It goes back to the work of Fujita \cite{fujita} that there is an exponent $p_{F}$ (critical Fujita exponent) for which the behavior of globally defined solutions of the reaction-diffusion model
\[
u_t-\Delta u=|u|^{p} \ \hspace{0.1cm}\mbox{in}\hspace{0.1cm}
    \mathbb{R}^{N}\times(0,\infty)
\]
are classified according to whether one is in the subcritical case ($1<p<p_{F}$), critical ($p=p_{F}$) or supercritical ($p_{F}<p<\infty)$. A very well-known argument allowing one to predict the correct value of $p_{F}$, whenever it is finite, for several semilinear parabolic equations consists of the study of the (linear) stability of the zero solution, see for instance \cite{Lev}. A more direct and simple approach pertaining to problems posed in the whole space and intrinsically connected to the scaling of the underlying equation has been proposed by the authors in \cite{Caz-Dick-Weiss}. More specifically, a scaling analysis on \eqref{eq:base} without the forcing term, that is,
\begin{equation}\label{eq:Hardy-Henon}
u_t+(-\Delta)^d u=|x|^{\alpha}|u|^{p} \ \hspace{0.1cm}\mbox{in}\hspace{0.1cm}
    \mathbb{R}^{N}\times(0,\infty),\hspace{0.2cm}u(0)=u_0
\end{equation}
shows that the only Lebesgue space which leaves invariant the norm of the appropriately rescaled initial data is $L^{p_c}(\RN)$, $p_{c}=\frac{N(p-1)}{2d+\alpha}$. Thus one may anticipate the critical Fujita exponent for \eqref{eq:Hardy-Henon} to be given by $1+\frac{2d+\alpha}{N}$ which corresponds to the value of $p$ for which $p_c=1$. The recent work \cite{G} asserts in turn that this observation is actually correct. Indeed, the global well-posedness theory for \eqref{eq:Hardy-Henon} among other questions has been investigated in the latter reference -- in particular, small data (in weak-Lebesgue space $L^{p_c,\infty}(\R^{N})$ global existence was obtained in $C_{b}([0,\infty);L^{p_c,\infty}(\R^{N}))$, $p_c=\frac{N(p-1)}{2d+\alpha}$ provided we require $\frac{N+\alpha}{N-2d}<p<\infty$ and $0<-\alpha<2d<N$ while under the conditions $u_0\in L^1(\RN)$ with $\displaystyle\int_{\R^N}u_0(x)dx>0$ and $1<p<1+\frac{2d+\alpha}{N}$, non-existence of weak solutions was established, see \cite{G}. Note that the previous well-posedness statement does not cover the range $(1+\frac{2d+\alpha}{N},\frac{N+\alpha}{N-2d})$, a restriction which essentially emanates from the choice of the functional setting together with the hypotheses needed for the applicability of certain key estimates heavily utilized in the analysis. This gap, however, could be closed by using the two-norms technique as illustrated in \cite{G}. Moving on, it is easy to see that the previous procedure fails to predict the critical exponent $p^{\star}$ for \eqref{eq:base} due to the presence of a forcing term. Nevertheless, as we shall see, granted the above knowledge on \eqref{eq:Hardy-Henon}, the behavior of the inhomogeneous term plays a crucial role in determining $p^{\star}$. Throughout this manuscript, unless otherwise specified we will assume that $\zeta : (0,\infty)\to (0,\infty)$ is a continuous function satisfying either
\begin{equation}
\label{eq:Z1}
\zeta(t)\,\underset{t\to 0}\sim\, t^{\sigma},\quad \zeta(t)\underset{t\to \infty}\sim\, t^m, \;\;\; (\sigma>-1,\,m\in\R),
\end{equation}
or
\begin{equation}
\label{eq:Z2}
\zeta(t)=t^{\sigma},\;\; \sigma>-1.
\end{equation}
 The study of the reaction-diffusion equation involving a source term has been considered by a number of authors. When the forcing term is purely space-dependent (i.e. $\zeta\equiv 1$) with $d=1$ and $\alpha=0$, it is known that the Fujita critical exponent reads
\[  p_{F}= \left\{
\begin{array}{ll}
      \infty\hspace{0.2cm}\mbox{if}\hspace{0.2cm} N=1,2 \\
      \frac{N}{N-2} \hspace{0.2cm}\mbox{if}\hspace{0.2cm} N>2
\end{array}
\right. \]
in the sense that for $p$ below $p_{F}$ and $\displaystyle\int_{\RN}\textbf{w}(x)dx>0$, Problem \eqref{eq:base} does not admit any global solution -- a conclusion which persists in the critical case $p=p_{F}$ under the additional assumptions on $\textbf{w}$ whilst for $p>p_{F}$ global-in-time solutions exist provided $\textbf{w}$ and $u_0$ have strong decay property at infinity. These statements can be found in \cite{BLZ} where various extensions and results related to other geometries were obtained, see also \cite{Zh2} for the geometrical setting given by a manifold. When the source term contains a polynomial in time just as considered in \eqref{eq:Z2} and for $d=1$, $\alpha=0$ it was proved in \cite{JKS} that the critical exponent is a function of $\sigma$ which for $N\geq 3$ has a discontinuity at zero,
\[  p_{F}(\sigma)= \left\{
\begin{array}{ll}
      \frac{N-2\sigma}{N-2\sigma-2}\hspace{0.2cm}\mbox{if}\hspace{0.2cm} -1<\sigma<0 \\
      \infty \hspace{0.2cm}\mbox{if}\hspace{0.2cm} \sigma>0.
\end{array}
\right. \]
The situation $d=1$ and $-2<\alpha<0$ has been investigated in \cite{Ma}. It was proved that for $N\geq 3$, solutions blow up in each of the following cases: $m\leq 0$, $p< \frac{N-2m+\alpha}{N-2m-2}$, $\int_{\RN}\textbf{w}dx>0$ and $p>1$. Furthermore, these results do not depend on the behavior of $\zeta$ at small times. However, the study of global-in-time solutions has not been taken into account and we wish to cover that in this work.

 Our results generalize earlier existing works and reveal new aspects of the solution. Prescribing initial data $u_0$ and $\textbf{w}$ small in Marcinkiewicz space gives rise to a global solution which lives in a small ball in the solution space and may as well be singular. As a byproduct, global existence holds in Lebesgue analogues of the previous functional setting whenever $u_0$ and $\textbf{w}$ are chosen small in the corresponding topologies. In particular, there are solutions which are initially singular and if additionally data are assumed to be radially symmetric (resp. radially monotone), then the solution enjoys the same properties. Moreover, when the principal operator in \eqref{eq:base} is the nonlocal fractional Laplacian ($d \in (0,1)$) or the Laplace operator $\Delta$, ($d=1$) then the solution is positive along with $u_0$ and $\textbf{w}$. In the process, we use an approach which in nature differs from those which have been employed in the aforementioned papers -- for instance, the use of Kato's two-norms technique which disallow unconditional uniqueness criterion and excludes a priori the existence of solutions which at time $t=0$ are singular -- our global wellposedness statement claims the converse. We essentially rely on smoothing effect type bounds in weak-Lebesgue spaces \cite{G} for the singular polyharmonic heat semigroup and its integral kernel combined with Yamazaki's type estimates \cite{Ya} which have a wide range of application -- they allow us among others to derive unconditional uniqueness of globally defined solutions for prescribed data in $L^{p_c}(\RN)$. Remarkably, we allow $u_0$ and $\mathbf{w}$ to be singular functions on $\RN$. Concerning blow-up results, our Theorem \ref{theo:non-existence} reinforces the fact that blow up whenever it occurs does not depend on the behavior of the forcing term at small times. To prove this theorem, we proceed via the test function method introduced by the authors in \cite{MP}.

We adopt the following notion of solution.
\begin{definition}\label{defn:weak-solution}
We say that $u$ is a global weak solution of \eqref{eq:base} if it satisfies the conditions
\begin{equation}\label{eq:weak-solution}
\nonumber u_0\in L^1_{loc}(\RN),\quad |x|^{\alpha}|u|^{p}\in L^1_{loc}(\RN\times (0,\infty))\end{equation} and
\begin{align*}
\int_{\RN\times (0,\infty)}u(-\partial_{t}\psi+(-\Delta)^d\psi)dx dt=\int_{\RN}u_0(x)\psi(x,0)dx+\int_{\RN\times (0,\infty)}|x|^{\alpha}|u|^p\psi dxdt+\\
\qquad{}\hspace{3cm}\int_{\RN\times (0,\infty)}\zeta(t)\textbf{w}\psi\,dxdt
\end{align*}
for all $\psi\in C^{\infty}_{0}(\RN\times [0,\infty))$.
\end{definition}

As is a standard practice, (\ref{eq:base}) is equivalent in an appropriate framework to the Duhamel formulation
\begin{equation}\label{eq:Duhamel}
u(x,t)=e^{-t(-\Delta)^{d}}u_0+
\int_{0}^{t}e^{-(t-s)(-\Delta)^{d}}
\left(|\cdot|^{\alpha}|u(s)|^p\right)\,ds
+\int_{0}^{t}\, \zeta(s)\,e^{-(t-s)(-\Delta)^{d}}{\mathbf w} \,ds,
\end{equation}
where $e^{-t(-\Delta)^{d}}$ is the linear semi-group generated by $(-\Delta)^{d}$. A solution to the integral equation \eqref{eq:Duhamel} is often called mild solution of \eqref{eq:base}. In our setting, one can see that a mild solution is a weak solution. See \cite{Terr1, Terr2, RT, PGL} for a discussion on the equivalence between the differential and integral formulation for the nonlinear heat equation.

In parallel to the questions discussed above, we also study local existence of solutions. Our first result in this direction deals with the case of initial data in the space of bounded continuous functions in $\RN$.
\begin{theorem}\label{theo:localwell-cont}
Suppose we have \eqref{eq:Z1}, ${\mathbf w}\in BC(\mathbb{R}^{N})$, $0<-\alpha<\min(2d,N)$ and $N\geq 2$. Given $u_0\in BC(\mathbb{R}^{N})$, there exists a time $T:=T(u_0)>0$ and a unique mild solution $u\in C([0,T];BC(\mathbb{R}^{N}))$ of Eq. \eqref{eq:base} such that $u(0)=u_0.$
\end{theorem}
Our next goal is to investigate the local well-posedness in Lebesgue spaces. The first result of this flavor was obtained in \cite{BTW} for \eqref{eq:base} with $d=1$ and without a forcing term.
\begin{theorem}\label{theo:localwell-Lp}
Let  $p_c=\frac{N(p-1)}{2d+\alpha}$. Assume $0<-\alpha<\min(2d, N)$ and $q>\max\big\{p_c,\frac{Np}{N+\alpha}\big\}$ with $q>1$. Suppose that \eqref{eq:Z1} holds, $u_0, {\mathbf w}\in L^q(\mathbb{R}^{N})$ and $N\geq 2$.Then, there exists $T:=T(u_0)>0$ and a unique function $u\in C([0,T]; L^{q}(\mathbb{R}^{N}))$ mild solution of Eq. \eqref{eq:base}. Let $T_{\max}$ be the maximal time of existence of this solution. If $T_{\max}<\infty$, then one has
$$\lim_{t\rightarrow T_{\max}}\|u(t)\|_{L^{q}(\mathbb{R}^{N})}=\infty.$$
\end{theorem}

\begin{remark}Observe that although the solution obtained in Theorem \ref{theo:localwell-cont} has the regularity of the initial data, it is not true in general that this solution enjoys a better regularity, see \cite{Wang} for this evidence in some special situation. A similar result can be stated for data in $C_{0}(\RN)$, the set of continuous function  vanishing at infinity in which case, the solution naturally belongs to $C_{0}(\RN)$. In theorem \ref{theo:localwell-Lp}, the fact that one can compare $p_c$ and $\frac{Np}{N+\alpha}$ leads to the question whether or not local solutions can be constructed if the condition imposed on $q$ fails to hold true. This involves three different cases namely $p_c>\frac{Np}{N+\alpha}$, $p_c=\frac{Np}{N+\alpha}$ and $p_c<\frac{Np}{N+\alpha}$ which are worth investigating. We do not pursue in this direction here. We point out, however, that local wellposedness in the critical case $p_c>\frac{Np}{N+\alpha}$ appears as a byproduct of Theorem \ref{theo:global existence} below.
\end{remark}
Next, motivated by the study of local-in-time solutions which are smooth $C^{\infty}$ instantly, we consider $\alpha>0$ and introduce the space $C_{\Lambda}(\mathbb{R}^{N})$ which collects continuous functions defined on $\mathbb{R}^{N}$ satisfying the condition
$\|\Lambda u\|_{L^{\infty}(\mathbb{R}^{N})}<\infty$ where $\Lambda(x)=(1+|x|)^{\frac{\alpha}{p -1}}$ and $p>1$. We endow this space with the norm  given by $$\|u\|_{\Lambda}=\|\Lambda u\|_{L^{\infty}(\mathbb{R}^{N})}.$$
The following result establishes local existence for initial data in $C_{\Lambda}(\mathbb{R}^{N})$.

\begin{theorem}\label{theo:localwell-decay} Let $\alpha> 0$ and $\zeta$ satisfying \eqref{eq:Z2}. Then, Eq. \eqref{eq:Duhamel} is locally wellposed. More precisely, given $u_0, {\mathbf w} \in C_{\Lambda}(\mathbb{R}^{N})$, there exist $T\in (0,\infty)$ and a unique solution $u$ of \eqref{eq:Duhamel} in $C([0,T]; C_{\Lambda}(\mathbb{R}^{N}))$. This solution in turn is classical and can be extended to a maximal interval $[0,T_{\max})$ where $T_{\max}\leq \infty$ and $\displaystyle \lim_{t\rightarrow T_{\max}}\|u(t)\|_{\Lambda}=\infty$ whenever $T_{\max}<\infty$.
\end{theorem}
Switching to the analysis of the global theory, our main results in this direction are the followings.
\begin{theorem}[Global existence]\label{theo:global existence}
Let $0<-\alpha<2d<N$. Granted \eqref{eq:Z2} with $-1<\sigma<0$, assume that $p\geq \frac{N-2d\sigma+\alpha}{N-2d\sigma-2d}$ and set $\ell=\frac{N p_c}{N+2(\sigma+1)d p_c}$.  Then for any $u_0\in L^{p_c,\infty}(\mathbb{R}^N)$ and ${\mathbf w}\in L^{\ell,\infty}(\mathbb{R}^N)$ with the property that $\|u_0\|_{L^{p_c,\infty}}+\|{\mathbf w}\|_{L^{\ell,\infty}}$ is sufficiently small,  Eq. \eqref{eq:Duhamel} admits a global-in-time solution $u$ which converges to $u_0$ as $t\rightarrow 0^{+}$ in $\mathcal{S}'(\RN)$. Moreover, this solution obeys the following properties.
\begin{enumerate}[label={($\textbf{a}_{\arabic*})$}]
\item\label{radial-symmetry}$(\textbf{Radial symmetry})$. If $u_0$ and ${\mathbf w}$ are radial functions, then the solution $u$ is radial in the variable $x$.
\item \label{monotonicity}$(\textbf{Radial monotonicity})$. Assume $u_0$ and ${\mathbf w}$ are radially nonincreasing. Then $u(t)$ is radially nonincreasing in $x$ for all $t\in [0,\infty)$.
\item\label{positivity} $(\textbf{Positiveness})$. Assume $d\in (0,1]$ and let $u_0$ and $\mathbf{w}$ be positive functions. Then the solution $u$ is also positive.
\end{enumerate}
\end{theorem}
\begin{remark}
The above theorem, in particular shows that the functions $u_0(x)=|x|^{-\frac{2d+\alpha}{p-1}}$ and $\mathbf{w}(x)=|x|^{-\frac{N}{\ell}}$ as well as their translated analogues $|x-x_0|^{-\frac{2d+\alpha}{p-1}}$ and $|x-x_0|^{-\frac{N}{\ell}}$ for each $x_0\in \RN$ are admissible choices as they belong to $L^{p_c,\infty}(\RN)$ and $L^{\ell,\infty}(\RN)$, respectively. A direct consequence of Theorem \ref{theo:global existence} is that for $\mathbf{w}\in L^{\ell}(\RN)$, there exist solutions which initially belong to the critical Lebesgue space $L^{p_c}(\RN)\subset L^{p_c,\infty}(\RN)$.
\end{remark}
\begin{theorem}[Unconditional uniqueness]\label{theo:unconditional-uniqueness}
Assume that $-\alpha\in (0,2d)$ and let $2d<N$. Further assume that $p>\frac{N+\alpha}{N-2d}$. If $u$ and $v$ are two solutions of \eqref{eq:Duhamel} in $C([0,\infty);L^{p_c}(\RN))$ obtained under smallness of $\|u_0\|_{L^{p_c}}+\|\mathbf{w}\|_{L^{\ell}},$ then $u=v$ a.e. in $\RN\times (0,\infty)$.
\end{theorem}
\begin{remark}
The result of Theorem \ref{theo:unconditional-uniqueness} is true for the inhomogeneous nonlinearity $|x|^\alpha\, |u|^{p-1}u$. Note that the same result was obtained in \cite{Tay} for $d=1$ and $\zeta\,{\mathbf w}=0$ and in \cite{G} for $d\in (0,1)\cup \N$ and $\zeta\,{\mathbf w}=0$. See \cite[Theorem 1.1 (ii), p. 4]{Tay} for $d=1$. The idea of the proof in \cite{G} is the same as in \cite{Tay}.
Here, we adapt the arguments used in \cite{Tay} in our setting. See also \cite[Proposition 3.7, p. 14]{Tay} and \cite[Lemma 4.2, p. 16]{Tay}.
\end{remark}
\begin{theorem}[Blow-up]\label{theo:non-existence} Suppose that $d\geq 1$ is an integer, $0<-\alpha<2d<N$, and that ${\mathbf w}$ belongs to $C_0(\mathbb{R}^N)\cap L^1(\mathbb{R}^N)$, obeys $\displaystyle\int_{\mathbb{R}^N}{\mathbf w}(x)dx>0$. Let $\zeta$ defined as in \eqref{eq:Z1}.
\begin{itemize}
\item[(i)]If $m\leq 0$ and $1<p<\frac{N-2dm+\alpha}{N-2m d-2d}$, then for any initial data $u_0\in C_0(\mathbb{R}^N)$, then Pb. \eqref{eq:base} has no global solution in the sense of definition \ref{defn:weak-solution}.
\item[(ii)]In case $m>0$, the same conclusion holds whenever $p>1$, i.e. Pb. \eqref{eq:base} has no global weak solution.
\end{itemize}
\end{theorem}
\begin{remark}
The statements of Theorem \ref{theo:non-existence} and Theorem \ref{theo:global existence} both reveals that \begin{equation*}
  p_F(\sigma) = \left\{
        \begin{array}{ll}
            \frac{N-2d\sigma+\alpha}{N-2\sigma d-2d} & \text{if}\hspace{0.1cm} -1 < \sigma<0 \\
            \infty & \text{if}\hspace{0.1cm}  \sigma>0
        \end{array}
    \right.
\end{equation*}
 separates the nonexistence/existence regime of global-in-time solutions of \eqref{eq:base} if $\zeta$ satisfies \eqref{eq:Z2}.
\end{remark}
\begin{remark}
The blow-up result does not cover the fractional case $d\in (0,1)$. The ideas of the proof when $d\in \mathbb{N}$ do not directly apply to the latter case because of the nature of the fractional Laplacian. See \eqref{eq:bound-on-d-Laplace} below. New ideas are therefore required.
\end{remark}

The rest of this paper is organized as follows. In the next section, we
recall some basic facts and useful tools. Section 3 is devoted to the proofs of our main results.
\section{Background materials}\label{prelim}
In this section we collect auxiliary results which will later find their usefulness and applicability in the proofs of our main results. Recall
the space $BC(\mathbb{R}^{N})$ which collects bounded and continuous functions in $\mathbb{R}^{N}$ and $C_0(\mathbb{R}^{N})$ denoting the space of continuous functions in $\mathbb{R}^{N}$ vanishing at infinity.

We will need some basic properties of Lorentz spaces. For $f:\RN\rightarrow \mathbb{R}$ a measurable function on $\RN$, denote by $f^{\ast}(\lambda)=\inf\{\tau>0:|\{|f|> \tau\}|\leq \lambda\}$ its decreasing function and let $f^{\ast\ast}(\tau)=\frac{1}{\tau}\int_{0}^{\tau}f^{\ast}(s)ds$. The Lorentz space $L^{p,q}(\RN)$, (weak-Lebesgue or Marcinkiewicz space in case $q=\infty$)  collects all measurable functions $f$ such that$\|f\|_{L^{p,q}}$ is finite,
\begin{equation*}
\|f\|_{L^{p,q}}=\left\{
\begin{array}{ll}
\bigg(\dfrac{q}{p}\displaystyle\int^{\infty}_{0}[t^{1/p}f^{\ast\ast}(t)]^q\dfrac{dt}{t}\bigg)^{1/q}, & \mbox{if}\hspace{0.2cm} p\in (1,\infty), \hspace{0.2cm}q\in [1,\infty) \\
\displaystyle\sup_{t>0}t^{1/p}f^{\ast\ast}(t), & \mbox{if} \hspace{0.2cm}p\in (1,\infty],\hspace{0.2cm} q=\infty.
\end{array}
\right.
\end{equation*}
These spaces increase with the second exponent and contain Lebesgue spaces as subspace. In fact, $L^{p,1}(\RN)\subset L^{p,q_1}(\RN)\subset L^{p}(\RN)\subset L^{p,q_2}(\RN)$, $1<p<\infty$, $1\leq q_1 \leq p\leq q_2\leq\infty$ with continuous injection. Moreover, $L^{p,q}(\RN)$ may alternatively be realized as the real interpolation space between the Lebesgue $L^1(\RN)$ and $L^{\infty}(\RN)$,
\[L^{p,q}(\RN)=[L^1(\RN),L^{\infty}(\RN)]_{1-1/p,q},\hspace{0.2cm}p\in (1,\infty), q\in [1,\infty]
\]
and more generally, for $1<p_1<p_2<\infty$, $\dfrac{1}{p}=\dfrac{1-\theta}{p_1}+\dfrac{\theta}{p_2}$, $\theta\in (0,1)$ and $1\leq q,q_1,q_2\leq \infty$ we have the interpolation identity
\[[L^{p_1,q_1}(\RN),L^{p_2,q_2}(\RN)]_{\theta,q}=L^{p,q}(\RN).
\]
Consider the linear polyharmonic heat equation
\begin{equation}\label{eq:polyharmoheat}
\partial_tu+(-\Delta)^{d}u=0 \hspace{0.2cm} \text{in} \hspace{0.2cm} \mathbb{R}^{N}\times (0,\infty).
\end{equation}
It is well known that the operator $(-\Delta)^d$, $d\geq 1$ is a generator of a semigroup $e^{-t(-\Delta)^d}$ whose kernel $E_d$ is smooth, radial and satisfies the scaling property
\begin{equation}\label{kernel}
E_{d}(x,t)=t^{-\frac{N}{2d}}E_d(t^{-\frac{1}{2d}}x,1).
\end{equation}
Hence, a solution of (\ref{eq:polyharmoheat}) subject to initial data $u(0)=u_0$ may be formally realized via convolution by \[u(x,t)=e^{-t(-\Delta)^{d}}u_0(x)=\left(E_{d}(\cdot,t)\ast u_0\right)(x)\] whenever this representation makes sense (e.g. when $u_0\in BC(\RN)$ or $u_0\in \mathscr{S}'(\RN)$ is a Schwartz distribution). Note that when $d\in (0,1)$, the operator $(-\Delta)^{d}$ defined via Fourier transform by $\widehat{(-\Delta)^{d}\varphi}(\xi)=|\xi|^{2d}\widehat{\varphi}(\xi)$ ($\varphi\in \mathcal{S}(\mathbb{R}^{N})$) also generates a semigroup whose kernel as an algebraic decay at infinity, see \cite{Jacob}. We recall the following $L^r-L^q$ estimate proved in \cite{ADE} by using the majorizing kernel established in \cite{GP}. See \cite[Proposition 6.1, p. 521]{ADE}.
\begin{proposition}\label{prop:smoothing-effect}
There exists a positive constant ${\mathcal H}_d$ such that for all $1\leq p\leq q\leq \infty$, we have
\begin{equation}\label{eq:smoothing-est}
\|e^{-t(-\Delta)^{d}}\varphi\|_{L^{q}}\leq {\mathcal H}_d\,t^{-\frac{N}{2d}(\frac{1}{p}-\frac{1}{q})}\|\varphi\|_{L^{p}}
\end{equation}
for all $t>0$ and $\varphi\in L^{p}(\mathbb{R}^{N})$.
\end{proposition}
Interpolating estimate \eqref{eq:smoothing-est} yields the following smoothing effect in Lorentz spaces
\begin{equation}\label{eq:smooth-effect-Lorentz}
\|e^{-t(-\Delta)^{d}}\varphi\|_{L^{q,r}}\leq {\mathcal H}_d\,t^{-\frac{N}{2d}(\frac{1}{p}-\frac{1}{q})}\|\varphi\|_{L^{p,r}}
\end{equation}
for all $t>0$, $\varphi\in L^{p,r}(\mathbb{R}^{N})$ and for $1\leq r\leq \infty$, $1<p\leq q<\infty$. For $N\geq 2$ and $\alpha\in (0,N)$, define the operator
\begin{equation*}
{\mathbf S}_{d,\alpha}(t)=e^{-t(-\Delta)^{d}}|\cdot|^{-\alpha}.
\end{equation*}

\begin{proposition}[\cite{G}]\label{prop:heat-hardy-kernel-est}
Let $N\geq 2$, $\alpha\in (0,N)$, $1<p_1,p_2\leq \infty$ and $q\in [1,\infty]$ such that
\begin{equation}
\label{eq:relation-p1-p2}
p_1>\frac{N}{N-\alpha},\hspace{0.5cm} p_2\geq \frac{N\,p_1}{N+\alpha\, p_1}.
\end{equation}
Then ${\mathbf S}_{d,\alpha}(t)$ has the following mapping properties:
\begin{enumerate}[label={($\textbf{A}_{\arabic*})$}]
\item \label{part a}${\mathbf S}_{d,\alpha}(t)$ maps continuously $L^{p_1}(\mathbb{R}^{n})$ into $C_0(\mathbb{R}^{n})$ for all $t>0$.
\item\label{part b} For any $\gamma \in \mathbb{N}^N_{0}$, ${\mathbf S}_{d,\alpha}(t)$ maps continuously $L^{p_1,\infty}(\mathbb{R}^{N})$ into $L^{p_2,q}(\mathbb{R}^{N})$ for all $t>0$. Moreover, there exits a positive constant $C:=C(N,p_1,p_2,d,\alpha,\gamma)$ such that
\begin{equation}\label{eq:mainest}
\|\partial^{\gamma}{\mathbf S}_{d,\alpha}(t)\varphi\|_{L^{p_2,q}}\leq Ct^{-\frac{N}{2d}(\frac{1}{p_1}-\frac{1}{p_2})-\frac{\alpha}{2d}-\frac{|\gamma|}{2d}}\|\varphi\|_{L^{p_1,\infty}}
\end{equation}
for all $t>0$ and $\varphi\in L^{p_1,\infty}(\RN)$.
\end{enumerate}
\end{proposition}
\begin{remark}
Proposition \ref{prop:heat-hardy-kernel-est} is known for $d=1, \gamma=0$ in \cite{BTW} for  Lebesgue spaces and in \cite[Proposition 3.3, p. 13]{Tay} for  Lorentz spaces. The original ideas of the proofs are like that of \cite{BTW, Tay}. See also \cite{G}.
\end{remark}

\begin{lemma}[\cite{G}]\label{lem:exp-est}
Given $\eta>0$, there exits a positive constant $C:=C(\eta,N)$ such that for $\beta\in [0,2]$ and $y\in \mathbb{R}^{N}$,
\begin{equation}\label{eq:exp-est}
\int_{\mathbb{R}^{N}}|g(x)|(1+|y-\beta x|)^{-\eta}dx\leq C(1+|y|)^{-\eta},
\end{equation}
where $g(x)=E_d(x,1)$ and $E_d$ is as in \eqref{kernel}.
\end{lemma}
Let $\phi:\RN\times (0,\infty)\rightarrow \mathbb{R}$. We introduce the nonlinear operator \[\mathscr{M}\phi(x)=\displaystyle\int_{0}^{\infty}{\mathbf S}_{d,\alpha}(s)|\phi(x,s)|^pds\] for a suitable $\phi$ such that the integral on the right hand side is meaningful.
\begin{lemma}[\cite{G}]\label{lem:Yamazaki's-bound}
Let $p>1$ and $0<\alpha<2d<N$. Given $k,q>1$, assume that \[k>\frac{Np}{N-\alpha}\hspace{0.2cm}\text{and}\hspace{0.2cm} \frac{1}{q}=\frac{p}{k}-\frac{2d-\alpha}{N}.\] There exists a constant $C>0$ such that
\begin{equation}\label{eq:1st-bound-on-M}
\|\mathscr{M}\phi\|_{L^{q,\infty}}\leq C\sup_{t>0}\|\phi(t)\|^{p}_{L^{k,\infty}}, \hspace{0.2cm} \phi\in L^{\infty}([0,\infty);L^{k,\infty}(\RN)).
\end{equation}
In particular we can take $k=\frac{N(p-1)}{2d-\alpha}$ in \eqref{eq:1st-bound-on-M} and if $\phi\in L^{\infty}([0,\infty);L^{p_c,\infty}\cap L^{r,\infty}(\RN))$ with $r>1$ and $1<r'<\frac{N}{2d-\alpha}$, $\frac{1}{r}+\frac{1}{r'}=1$, then we have
\begin{equation}\label{eq:snd-bound-on-M}
\|\mathscr{M}\phi\|_{L^{r,\infty}}\leq\, C\sup_{t>0}\|\phi(t)\|_{L^{r,\infty}}\sup_{t>0}\|\phi(t)\|^{p-1}_{L^{p_c,\infty}}.
\end{equation}
\end{lemma}
\section{Proofs of the main results}
Starting with the proofs of theorems pertaining to the local theory, recall that $\mathbf{S}_{d,\alpha}(t)$ is a $C_0$ semigroup on Lebesgue spaces and on $BC(\RN)$ so that we can systematically appeal to the abstract result in \cite{Weissler}.  For this end, we define for $0<-\alpha<N$ the map \[\mathscr{K}_tu=\mathbf{S}_{d,-\alpha}(t)|u|^p+\mathbf{S}_{d,0}(t)\zeta(\cdot)\mathbf{w}, \hspace{0.2cm} t>0.\]
\subsection{Proof of Theorem \ref{theo:localwell-cont}}
Given $a,b\in \mathbb{R}$, we can find a constant $\eta:=\eta(p)>0$ such that
\begin{equation}
\label{eq:nonlinearity-pointwise-ineq}
\big||a|^p-|b|^p\big|\leq \eta |a-b|(|a|^{p-1}+ |b|^{p-1}),\hspace{0.2cm}p>1.
\end{equation}
By invoking Proposition \ref{prop:heat-hardy-kernel-est} with $p_2=p_1=\infty$, it is easy to verify that $\mathscr{K}_{t}$ is bounded on the space $BC(\RN)$ and the following bound holds, namely $\|\mathbf{S}_{d,0}u_0\|_{L^\infty}\leq c\|u_0\|_{L^{\infty}}$. Next, call $\overline{B}_{L}(0)\subset BC(\RN)$ the ball with center at zero and radius $L>0$, we prove that $\mathscr{K}_t$ is Lipschitz continuous on $\overline{B}_{L}(0)$ for each $t>0$. Take $u, v\in \overline{B}_{L}(0)$, by \eqref{eq:nonlinearity-pointwise-ineq} and making use of Proposition \ref{prop:heat-hardy-kernel-est} once again we arrive at
\begin{align*}
\|\mathscr{K}_tu-\mathscr{K}_tv\|_{L^{\infty}}\leq c_{L}(t)\|u-v\|_{L^{\infty}}
\end{align*}
with the Lipschitz constant given by $c_{L}(t)=CL^{p-1}t^{\frac{\alpha}{2d}}\in L^{1}(0,\varepsilon)$ for some $\varepsilon>0$ and for each $L>0$ since $-\alpha<2d$. Moreover, the kernel $E_{d}$ has the property
\[E_{d}(x,t)=\int_{\RN}E_{d}(x-y,t-s)E_{d}(y,s)dy
\]
for all $x\in \RN$ and $s,t>0$; which obviously implies $\mathbf{S}_{d,0}(t)\phi=\mathbf{S}_{d,0}(t-s)[\mathbf{S}_{d,0}(s)\phi]$, $0<s<t$. Hence,
\begin{align*}
\mathbf{S}_{d,0}(s)\mathcal{K}_{t}u&=\mathbf{S}_{d,0}(s)\mathbf{S}_{d,0}(t)|\cdot|^{\alpha}|u|^p+\mathbf{S}_{d,0}(s)\mathbf{S}_{d,0}(t)\zeta(\cdot)\mathbf{w}\\
&=\mathbf{S}_{d,0}(t+s)|\cdot|^{\alpha}|u|^p+\mathbf{S}_{d,0}(t+s)\zeta(\cdot)\mathbf{w}\\
&=\mathbf{S}_{d,-\alpha}(t+s)|u|^p+\mathbf{S}_{d,0}(t+s)\zeta(\cdot)\mathbf{w}\\
&=\mathcal{K}_{t+s}u.
\end{align*}
On the other hand, we have $\mathcal{K}_{t}(0)=\mathbf{S}_{d,0}(t)\zeta(\cdot)\mathbf{w}$, so that for some  $\varepsilon>0$ and by using \eqref{eq:Z1}, there holds
\begin{align*}\big\|\|\mathscr{K}_t(0)\|_{L^{\infty}}\big\|_{L^1(0,\varepsilon)}&= \big\|\|\mathbf{S}_{d,0}(t)\zeta(\cdot)\mathbf{w}\|_{L^{\infty}}\big\|_{L^1(0,\varepsilon)}\\
&\leq C\|\mathbf{w}\|_{L^{\infty}}\|\zeta\|_{L^1(0,\varepsilon)}\leq C.
\end{align*}
At this point, a simple application of Theorem 1 in \cite{Weissler} yields the desired conclusion.


\subsection{Proof of Theorem \ref{theo:localwell-Lp}}
Here, the reasoning goes in the same spirit as in the above proof. Let $u_0,\mathbf{w}\in L^q(\RN)$. Then granted  \eqref{eq:Z1}, under the condition $q>\frac{Np}{2d+\alpha}$, we can apply Proposition \ref{prop:heat-hardy-kernel-est} with $(p_1,p_2)=(q/p,q)$ to obtain that $\mathscr{K}_t$ is bounded on $L^{q}(\RN)$ for each $t>0$ and there holds $\|\mathbf{S}_{d,0}u_0\|_{L^q}\leq C\|u_0\|_{L^q}$ by Proposition \ref{prop:smoothing-effect}. We eventually want to verify that the following conditions are fulfilled:
\begin{enumerate}
    \item $\mathscr{K}_t:L^q(\RN)\rightarrow L^q(\RN)$ is a locally Lipschitz map with constant $c_L(t)\in L^1(0,\varepsilon)$ for each $L>0$ and for some $\varepsilon>0$.
    \item The map $t\mapsto \|\mathscr{K}_t(0)\|_{L^{q}}$ belongs to $L^1(0,\varepsilon)$ for some $\varepsilon>0.$
    \item $\mathbf{S}_{d,0}(s)\mathscr{K}_t=\mathscr{K}_{t+s}$ for $s,t>0.$
\end{enumerate}
The third property has already been established in the previous lines. As regards the second, for an arbitrary $\varepsilon>0$, we have that
\begin{align*}\big\|\|\mathscr{K}_t(0)\|_{L^{q}}\big\|_{L^1(0,\varepsilon)}&= \big\|\|\mathbf{S}_{d,0}(t)\zeta(\cdot)\mathbf{w}\|_{L^{q}}\big\|_{L^1(0,\varepsilon)}\\
&\leq C\|\mathbf{w}\|_{L^{q}}\|\zeta\|_{L^1(0,\varepsilon)}\leq C.
\end{align*}
Let us now show that the first condition is equally satisfied. Fix $L>0$ and denote by $\overline{B}_L(0)$ the closed ball in $L^q(\RN)$ centered at the origin and with radius $L$. We compute
\begin{align*}
\|\mathscr{K}_tu-\mathscr{K}_tv\|_{L^{q}}&=\|\mathbf{S}_{d,-\alpha}(t)(|u|^p-|v|^p)\|_{L^{q}}\\
&\leq c_{L}(t)\|u-v\|_{L^{q}}
\end{align*}
with $c_L(t)=CL^{p-1}t^{-\frac{N(p-1)}{2qd}+\frac{\alpha}{2d}}$. Note that in order to get this bound, we use the pointwise estimate \eqref{eq:nonlinearity-pointwise-ineq} and invoked Proposition \ref{prop:heat-hardy-kernel-est} with $(p_1,p_2)=(\frac{q}{p},q)$. Since $q>\frac{N(p-1)}{2d+\alpha}=p_c$, it follows that $c_L(t)\in L^1(0,\varepsilon)$ for all $L>0$ and for all $\varepsilon>0$. We are now in position to apply Theorem 1 in \cite{Weissler} which yields the statement of Theorem \ref{theo:localwell-Lp}.

\subsection{Proof of Theorem \ref{theo:localwell-decay}}
We turn to the proof of the local well-posedness result subject to data in $C_{\Lambda}(\mathbb{R}^{n})$. Recall $\Lambda(x)=(1+|x|)^{\frac{\alpha}{p-1}}$. In this case, we proceed in a slightly different fashion with the aim of reaching the hypotheses that will allow an application of a Banach fixed point theorem.

Let $0<T<1$ be a time to be chosen later and $\delta>0$. Consider the space \[
X_{T,\Lambda}=\bigg\{u\in C([0,T];C_{\Lambda}(\mathbb{R}^{N})):\|\Lambda u(t)\|_{L^{\infty}(\mathbb{R}^{N})}\leq \delta, t\in (0,T)\bigg\}\]
on which we put the metric $d(u,v)=\displaystyle\sup_{t\in [0,T]}\|\Lambda (u(t)-v(t))\|_{L^{\infty}(\mathbb{R}^{N})}$. We will prove in the sequel that the operator $\mathcal{Q}$ defined as
\begin{equation}\label{eq:integral-eq}
\mathcal{Q}u=\mathbf{S}_{d,0}(t)u_0+\int_{0}^t\mathbf{S}_{d,-\alpha}(t-s)|u(s)|^pds+\int_{0}^t\mathbf{S}_{d,0}(t-s)\zeta(s)\mathbf{w}ds
\end{equation}
possesses a fixed point in $C([0,T];C_{\Lambda}(\mathbb{R}^{N}))$. We first show that $\mathcal{Q}$ is a self-mapping from $X_{T,\Lambda}$ onto itself. Using Lemma \ref{lem:exp-est}, we have
\begin{align*}
|\mathbf{S}_{d,0}(t)u_0|&=|e^{-t(-\Delta)^{d}}u_0|=\bigg|\int_{\mathbb{R}^{N}}\,E_d(x-y,t)u_0(y)\,dy\bigg|\\
&\leq \int_{\mathbb{R}^{N}}\,|E_d(x-y,t)||u_0(y)|dy\\
&\leq \int_{\mathbb{R}^{N}}\,t^{\frac{-N}{2d}}|g(t^{-\frac{1}{2d}}(x-y))||u_0(y)|dy\\
&\leq \|u_0\|_{C_{\Lambda}(\mathbb{R}^{N})}\,\int_{\mathbb{R}^{N}}\,t^{\frac{-N}{2d}}|g(t^{-\frac{1}{2d}}(x-y))|
\Lambda^{-1}(y)\,dy\\
&\leq \|u_0\|_{C_{\Lambda}(\mathbb{R}^{N})}\,\int_{\mathbb{R}^{N}}|g(y)|\left(1+|x-t^{1/2d}\,y|\right)
^{-\frac{\alpha}{p-1}}\,dy\\
&\leq C(\alpha,p,N)\|u_0\|_{C_{\Lambda}(\mathbb{R}^{N})}\,(1+|x|)^{-\frac{\alpha}{p-1}}.
\end{align*}
Similarly, setting $Bu(t)=\displaystyle\int_{0}^{t}e^{-(t-s)(-\Delta)^{d}}|\cdot|^{\alpha}|u(s)|^{p}\,ds$, $u\in X_{T,\Lambda}$, an application of Lemma \ref{lem:exp-est} once again permits us to write
 \begin{align*}
 |B(t)u|&=\bigg|\int_{0}^{t}\int_{\mathbb{R}^{N}}\,(t-s)^{-\frac{N}{2d}}g\big((t-s)^{-\frac{1}{2d}}(x-y)\big)
 |y|^{\alpha}|u(s)|^{p}\,dyds\bigg|\\
 &\leq C\int_{0}^{t}\int_{\mathbb{R}^{N}}\,(t-s)^{-\frac{N}{2d}}|g\big((t-s)^{-\frac{1}{2d}}(x-y)\big)||y|^{\alpha}|\Lambda(y)u(y,s)|^{p}\Lambda^{-p}(y)\,dyds\\
&\leq C \delta^{p} \int_{0}^{t}\,(t-s)^{-\frac{N}{2d}}
\int_{\mathbb{R}^{N}}\,|g\big((t-s)^{-\frac{1}{2d}}(x-y)\big)|y|^{\alpha}(1+|y|)^{-\frac{\alpha\,p}{p-1}}\,dyds\\
&\leq C\delta^{p} \int_{0}^{t}\int_{\mathbb{R}^{N}}\,|g(z)|(1+|x-(t-s)^{\frac{1}{2d}}z|)^{-\frac{\alpha}{p-1}}\,dzds\\
&\leq C\delta^{p}\int_{0}^{T}\,(1+|x|)^{-\frac{\alpha}{p-1}}\,ds\\
&\leq C\delta^{p}(1+|x|)^{-\frac{\alpha}{p-1}}T.
\end{align*}
Finally, setting $Du(t)=\displaystyle\int_{0}^{t}\,e^{-(t-s)(-\Delta)^{d}}{\mathbf w}\zeta(s)\,ds$, we invoke Lemma \ref{lem:exp-est} to arrive at
  \begin{align*} |D(t)u|&=\bigg| \int_{0}^{t}\int_{\mathbb{R}^{N}}\,(t-s)^{-\frac{N}{2d}}g\big((t-s)^{-\frac{1}{2d}}(x-y)\big)|{\mathbf w}(y)|\zeta(s)\,dyds\bigg|\\
  &\leq \|{\mathbf w}\|_{C_{\Lambda}(\mathbb{R}^{N})}\int_{0}^{t}\int_{\mathbb{R}^{N}}\,
  (t-s)^{-\frac{N}{2d}}\big|g\big((t-s)^{-\frac{1}{2d}}(x-y)\big)\big||\Lambda(y)|^{-1}\zeta(s)\,dyds\\
 &\leq  \|{\mathbf w}\|_{C_{\Lambda}(\mathbb{R}^{N})}\int_{0}^{t}\int_{\mathbb{R}^{N}}\,|g(z)|\zeta(s)
 (1+|x-(t-s)^{\frac{1}{2d}}z|^{-\frac{\alpha}{p-1}}\,dzds\\
 &\leq C\|{\mathbf w}\|_{C_{\Lambda}(\mathbb{R}^{N})}\frac{T^{\sigma+1}}{\sigma+1}\Lambda^{-1}(x)
\end{align*}
so that $\|\Lambda D u(t)\|_{L^\infty}\leq C\|{\mathbf w}\|_{C_{\Lambda}(\mathbb{R}^{N})}\dfrac{T^{\sigma+1}}{\sigma+1}$. Gluing together all the above estimates, we find that \begin{equation}\label{eq:glued-eq}\|\Lambda \mathcal{Q}u(t)\|_{L^\infty}\leq C\big(\|u_0\|_{C_{\Lambda}(\mathbb{R}^{N})}+\delta^{p}T+\frac{T^{\sigma+1}}{\sigma+1}\|{\mathbf w}\|_{C_{\Lambda}(\mathbb{R}^{N})}\big),\;\; 0\leq t\leq T,
\end{equation}
where $C$ is a constant depending on $N,\alpha$ and $p$. Let $\delta>0$ such that $C\|u_0\|_{C_{\Lambda}(\mathbb{R}^{N})}<\delta$. Pick $T<1$ with the property that $C\big(\|u_0\|_{C_{\Lambda}(\mathbb{R}^{N})}+\delta^{p}T+\frac{T^{\sigma+1}}{\sigma+1}\|{\mathbf w}\|_{C_{\Lambda}(\mathbb{R}^{N})}\big)\leq \delta$ to draw the conclusion that $\mathcal{Q}$ maps $X_{T,\Lambda}$ into itself. Next, let $u,v\in X_{T,\Lambda}$ and let's additionally impose the condition $C\delta^{p-1}T<1$ where $C$ is that constant appearing in \eqref{eq:glued-eq}. Then $\mathcal{Q}$ is a contraction. To see this write
\begin{align*}
|\mathcal{Q}u(t)-\mathcal{Q}v(t)|&=\bigg|\int_{0}^{t}\int_{\mathbb{R}^{N}}E_{d}(x-y,t-s)|y|^{\alpha}
\big||u(s)|^{p}-|v(s)|^{p}\big|\,dyds\bigg|\\
&\leq C\int_{0}^{t}\int_{\mathbb{R}^{N}}|E_{d}(x-y,t-s)||y|^{\alpha}|u-v|(|u|^{p-1}+|v|^{p-1})\,dyds\\
&\leq C\bigg(\int_{0}^{t}\int_{\mathbb{R}^{N}}|E_{d}(x-y,t-s)||y|^{\alpha}|u-v|\Lambda(y)|\Lambda u|^{p-1}\Lambda^{-p}(y)dyds+\\
&\qquad \int_{0}^{t}\int_{\mathbb{R}^{N}}|E_{d}(x-y,t-s)||y|^{\alpha}|u-v|\Lambda(y)|\Lambda v|^{p-1}\Lambda^{-p}(y)\,dyds\bigg)\\
&\leq Cd(u,v)\bigg(\|\Lambda u(t)\|^{p-1}_{L^\infty}\int_{0}^{t}\int_{\mathbb{R}^{N}}|g(z)|
|x-(t-s)^{\frac{1}{2d}}z|^{-1}dzds+\\
&\qquad\qquad\qquad\qquad \|\Lambda v(t)\|^{p-1}_{L^\infty}\int_{0}^{t}\int_{\mathbb{R}^{N}}|g(z)||x-(t-s)^{\frac{1}{2d}}z|^{-1}\,dzds\bigg).
\end{align*}
This yields, utilizing Lemma \ref{lem:exp-est} the bound
\[
d\left(\mathcal{Q}u, \mathcal{Q}v\right)\leq C\delta^{p-1}T d(u,v).
\]
Existence of a solution as claimed in Theorem \ref{theo:localwell-decay} is now a simple consequence of the Banach fixed point theorem. The blow up alternative is verified using standard arguments. Since $C_{\Lambda}(\RN)\subset C_{b}(\RN)$, we have that $u\in L^{\infty}([0,T];C_{b}(\RN))$ and by parabolic regularity theory, $u$ satisfies \eqref{eq:base} in the classical sense. The proof of Theorem \ref{theo:localwell-decay} is now complete.

Our next focus is on the proofs of global existence and non-existence results. Recall the definition of the beta function
\[\mathcal{B}(a,b)=\int_{0}^{1}s^{a-1}(1-s)^{b-1}ds=\frac{\Gamma(a)\Gamma(b)}{\Gamma(a+b)},\quad a,b>0.\]
where $\Gamma$ is the standard gamma function.

\subsection{Proof of Theorem \ref{theo:global existence}}
We will distinguish between two cases assuming first that $p>\frac{N+\alpha}{N-2d}$ with $-\alpha<2d<N$. The second case $p\in [p_{F},\frac{N+\alpha}{N-2d}]$ will be handle subsequently. Let $u_0\in L^{p_c,\infty}(\RN)$, $\mathbf{w}\in L^{\ell,\infty}(\RN)$ such that $\|u_0\|_{L^{p_c,\infty}}+\|\mathbf{w}\|_{L^{\ell,\infty}}<\varepsilon_0$ for some $\varepsilon_0>0$. We wish to show that the equation
\begin{equation}\label{eq:fixed-pt-eq}
u=w+\mathcal{F}(u) \hspace{0.2cm}\mbox{in}\hspace{0.2cm}\RN\times (0,\infty); \quad w=\mathbf{S}_{d,0}(t)u_0+\int_0^t\mathbf{S}_{d,0}(t-s)\zeta(s)\mathbf{w} ds
\end{equation}
has a unique fixed point in a closed ball $\overline{B}_{2\varepsilon}(0)$ of $C_{b}((0,\infty);L^{p_c,\infty})$, $\varepsilon=\varepsilon(\varepsilon_0)$ where for $t>0$, \[\mathcal{F}(u)(t)=\int_{0}^{t}\mathbf{S}_{d,-\alpha}(t-s)|u(s)|^p ds.
\]
Apply Lemma \ref{lem:Yamazaki's-bound} with $k=p_c$ (i.e $q=p_c$) in \eqref{eq:1st-bound-on-M} and
$|\phi(s)|^p=\begin{cases}
|u(t-s)|^p\hspace{0.2cm}\text{if} \hspace{0.2cm}s\in (0,t)\\
0 \hspace{0.2cm}\text{if} \hspace{0.2cm}s\geq t
\end{cases}$ to obtain the estimate
\begin{equation}\label{eq:self-mapping-bound}\|\mathcal{F}u\|_{L^{\infty}((0,\infty);L^{p_c,\infty})}\leq C\sup_{t>0}\|u(t)\|^{p}_{L^{p_c,\infty}}.
\end{equation}
Arguing in a similar fashion, we have that for $u$ and $v$ in $L^{\infty}((0,\infty);L^{p_c,\infty})$,
\begin{equation}\label{eq:contraction-bound}
\|\mathcal{F}(u)-\mathcal{F}(v)\|_{L^{\infty}(\mathbb{R}_{+};L^{p_c,\infty})}\leq C\|u-v\|_{L^{\infty}(\mathbb{R}_{+};L^{p_c,\infty})}(\|u\|^{p-1}_{L^{\infty}(\mathbb{R}_{+};L^{p_c,\infty})}+\|v\|^{p-1}_{L^{\infty}((0,\infty);L^{p_c,\infty})}).
\end{equation}
Observe that $\ell=\frac{Np_c}{N+2(\sigma+1)dp_c}>1$ in view of the condition imposed on $p$ and that
\begin{align*}
\|w\|_{L^{p_c,\infty}}\leq C\|u_0\|_{L^{p_c,\infty}}+\bigg\|\int_0^t\mathbf{S}_{d,0}(t-s)\zeta(s)\mathbf{w}ds\bigg\|_{L^{p_c,\infty}}
\end{align*}
where the second term is further estimated via duality as follows
\begin{align*}
\bigg\|\int_0^t\mathbf{S}_{d,0}(t-s)\zeta(s)\mathbf{w}ds\bigg\|_{L^{p_c,\infty}}&=\sup_{\substack{\psi\in L^{p_{c}',1}\\
	\|\psi\|_{L^{p_{c}',1}}=1}}\bigg|\langle\int_0^t\mathbf{S}_{d,0}(t-s)\zeta(s)\mathbf{w}ds,\psi\rangle\bigg|\\
&\leq \sup_{\substack{\psi\in L^{p_{c}',1}\\
		\|\psi\|_{L^{p_{c}',1}}=1}}\bigg| \int_{0}^{t}\int_{\RN}\mathbf{S}_{d,0}(t-s)\zeta(s)\mathbf{w}\psi(x)dxds\bigg|\\
&\leq \sup_{\substack{\psi\in L^{p_{c}',1}\\
		\|\psi\|_{L^{p_{c}',1}}=1}}\bigg(\int_{0}^{t}\|\mathbf{S}_{d,0}(t-s)\psi\|_{L^{\ell',1}}\|\zeta(s)\mathbf{w}\|_{L^{\ell,\infty}}ds\bigg)\\
&\leq C\sup_{\substack{\psi\in L^{p_{c}',1}\\
		 \|\psi\|_{L^{p_{c}',1}}=1}}\bigg(\int_{0}^{t}(t-s)^{-\frac{N}{2d}(\frac{1}{p_{c}'}-\frac{1}{\ell'})}\|\psi\|_{L^{p_{c}',1}}\|\zeta(s)\mathbf{w}\|_{L^{\ell,\infty}}ds\bigg)\\
&\leq C\sup_{\substack{\psi\in L^{p_{c}',1}\\
		 \|\psi\|_{L^{p_{c}',1}}=1}}\|\psi\|_{L^{p_{c}',1}}\|\mathbf{w}\|_{L^{\ell,\infty}}\int_{0}^{t}(t-s)^{-\frac{N}{2d}(\frac{1}{\ell}-\frac{1}{p_{c}})}\zeta(s)ds\\
&\leq C\sup_{\substack{\psi\in L^{p_{c}',1}\\
		\|\psi\|_{L^{p_{c}',1}}=1}}\|\psi\|_{L^{p_{c}',1}}\|\mathbf{w}\|_{L^{\ell,\infty}}\int_{0}^{1}\tau^{\sigma}(1-\tau)^{-(\sigma+1)}d\tau\\
&\leq C\|\mathbf{w}\|_{L^{\ell,\infty}}\mathcal{B}(\sigma+1,-\sigma)		
\end{align*}
where we have systematically used Fubini's Theorem, smoothing effect in Lorentz space \eqref{eq:smooth-effect-Lorentz}, the generalized H\"{o}lder's inequality and in the estimate before the last we performed the change of variables $\tau=s/t$ together with the fact that $\frac{1}{p_c'}-\frac{1}{\ell'}=\frac{2(\sigma+1)d}{N}$, $-1<\sigma<0$ and $\ell'$ being the conjugate exponent of $\ell$. Hence, we have that \[\|w\|_{L^{p_c,\infty}}\leq C_1(\|u_0\|_{L^{p_c,\infty}}+\|\mathbf{w}\|_{L^{\ell,\infty}})\leq C_1\varepsilon_0.\]
Put $\varepsilon=C_1\varepsilon_0$ and $R=(2^p\,C)^{-\frac{1}{p-1}}$. Assuming that $\varepsilon<R$, we deduce from \eqref{eq:self-mapping-bound} and \eqref{eq:contraction-bound} that Equation \eqref{eq:fixed-pt-eq} has a unique fixed point in $\overline{B}_{2\varepsilon}(0)$ which indeed is the  solution of \eqref{eq:Duhamel} we were looking for.
Now assume that $p_{F}\leq p\leq \frac{N+\alpha}{N-2d}$. We shall argue differently in this case since the strategy employed before will clearly fail to work well here. This being said, pick a number $r>1$ such that the inequality below is satisfied,
\begin{equation}\label{eq:cond-on-exponents}
\max\bigg\{\frac{\alpha p+2d}{Np(p-1)},\frac{1}{p_c}+\frac{2d\sigma}{N}\bigg\}<\frac{1}{r}<\min\bigg\{\frac{1}{p_c},\frac{N+\alpha}{Np}\bigg\},\quad r>p.
\end{equation}
This choice is possible. In fact, in view of the condition $p\geq p_F$, one can easily check that all inequalities in \eqref{eq:cond-on-exponents} are satisfied expect possibly for $\dfrac{1}{p_c}+\dfrac{2d\sigma}{N}<\dfrac{N+\alpha}{Np}$ which we verify as follows. The function $f(a)=2d\sigma p^2+p(2d-2d\sigma -a)+a+\alpha$ defined for $a\geq \rho=2d\sigma+\dfrac{2dp+\alpha}{p-1}$  is decreasing so that in particular, $f(N)\leq f(\rho)<0$ because $p\geq p_{F}=\frac{N-2d\sigma+\alpha}{N-2d(\sigma+1)}$.
Moving on, observe that from \eqref{eq:cond-on-exponents} follows the restrictions $1\leq \ell<p_c< r$ and if we put
\[\mu=\frac{N}{2d}\bigg(\frac{1}{p_c}-\frac{1}{r}\bigg),
\]
then  $0<\mu<\dfrac{1}{p}$. Moreover, we have
\[\mu=\frac{N(p-1)}{2rd}-\frac{\alpha}{2d}+p\mu-1=\frac{N}{2d}\bigg(\frac{1}{\ell}-\frac{1}{r}\bigg)-\sigma-1.
\]Introduce the function space $\mathbf{X}$ defined by
\[\mathbf{X}=\bigg\{u\in C_{b}\big((0,\infty);L^{p_c,\infty}(\RN)\big):t^{\mu}u\in C_{b}\big((0,\infty);L^{r,\infty}(\RN)\bigg\}
\]
which carries out the structure of a complete metric space when equipped with the distance $d(u,v)=\displaystyle\sup_{t>0}t^{\mu}\|u(t)-v(t)\|_{L^{r,\infty}}:=\|u-v\|_{\mathbf{X}}$. We want to show that  Eq. \eqref{eq:fixed-pt-eq} has a fixed point in $\mathbf{X}$ which is unique in $\overline{B}_{K}(0)\subset \mathbf{X}$ a closed ball centered at the origin and with radius $K>0$ sufficiently small.
We estimate separately each of the terms of the right hand side of \eqref{eq:fixed-pt-eq}. Making use of the smoothing estimate \eqref{eq:smooth-effect-Lorentz} we obtain
\begin{equation*}
\|\mathbf{S}_{d,0}(t)u_0\|_{L^{r,\infty}}\leq Ct^{-\frac{N}{2d}\big(\frac{1}{p_c}-\frac{1}{r}\big)}\|u_0\|_{L^{p_c,\infty}}\leq Ct^{-\mu}\|u_0\|_{L^{p_c,\infty}}.
\end{equation*}
Since $r>\frac{Np}{N+\alpha}$, it follows from Proposition \ref{prop:heat-hardy-kernel-est} that
\begin{align*}
\bigg\|\int_{0}^{t}\mathbf{S}_{d,-\alpha}(t-s)|u(s)|^{p}ds\bigg\|_{L^{r,\infty}}&\leq C\int_{0}^{t}
(t-s)^{-\frac{N}{2d}\big(\frac{p}{r}-\frac{1}{r}\big)+\frac{\alpha}{2d}}\||u|^{p}\|_{L^{r/p,\infty}}ds\\
&\leq C\int_{0}^{t}\,(t-s)^{-\frac{N}{2d}\big(\frac{p}{r}-\frac{1}{r}\big)+\frac{\alpha}{2d}}s^{-p\mu}(s^{\mu}\|u\|_{L^{r,\infty}})^{p}ds\\
&\leq C(\sup_{t>0}t^{\mu}\|u\|_{L^{r,\infty}})^{p}\int_{0}^{t}\,(t-s)^{-\frac{N(p-1)}{2rd}+\frac{\alpha}{2d}}s^{-p\mu}ds\\
&\leq Ct^{-\frac{N(p-1)}{2rd}+\frac{\alpha}{2d}-p\mu+1}\|u\|^{p}_{\mathbf{X}}\int_{0}^{1}(1-s)^{-\frac{N(p-1)}{2rd}+\frac{\alpha}{2d}}s^{-p\mu}ds\\
&\leq Ct^{-\mu}\|u\|^{p}_{\mathbf{X}}\mathcal{B}\bigg(1-p\mu,1-\frac{N(p-1)}{2rd}+\frac{\alpha}{2d}\bigg)\\
&\leq Ct^{-\mu}\|u\|^{p}_{\mathbf{X}}.
\end{align*}
To estimate the last term, we utilize \eqref{eq:smooth-effect-Lorentz} and proceed as follows
\begin{align*}
\bigg\|\int_{0}^{t}\mathbf{S}_{d,0}(t-s)\mathbf{w}\zeta(s)\bigg\|_{L^{r,\infty}}ds&\leq C\int_{0}^{t}(t-s)^{-\frac{N}{2d}\big(\frac{1}{\ell}-\frac{1}{r}\big)}\|\mathbf{w}\|_{L^{\ell,\infty}}
\zeta(s)ds\\
&\leq C\|\mathbf{w}\|_{L^{\ell,\infty}}t^{-\frac{N}{2d}\big(\frac{1}{\ell}-\frac{1}{r}\big)+\sigma+1}
\int_{0}^{1}(1-s)^{-\frac{N}{2d}\big(\frac{1}{\ell}-\frac{1}{r}\big)}s^{\sigma}ds\\
&\leq C\|\mathbf{w}\|_{L^{\ell,\infty}}t^{-\frac{N}{2d}\big(\frac{1}{\ell}-\frac{1}{r}\big)+\sigma+1}
\mathcal{B}\bigg(\sigma+1,1-\frac{N}{2d}\big(\frac{1}{\ell}-\frac{1}{r}\big)\bigg)\\
&\leq C\|\mathbf{w}\|_{L^{\ell,\infty}}t^{-\mu}.
\end{align*}
Summarizing, if $\mathcal{Q}u=w+\mathcal{F}(u), \hspace{0.1cm} w=\mathbf{S}_{d,0}(t)u_0+\int_0^t\mathbf{S}_{d,0}(t-s)\zeta(s)\mathbf{w}$, then  \[\displaystyle
\sup_{t>0}t^{\mu}\|\mathcal{Q}u(t)\|_{L^{r,\infty}}\leq C\left(\|u_0\|_{L^{p_c,\infty}}+K^{p}+\|\mathbf{w}\|_{L^{\ell,\infty}}\right)\leq C\varepsilon_0.
\]
 Upon taking $\varepsilon_0$ and $K>0$ sufficiently small, one can achieve $C(\|u_0\|_{L^{p_c}}+\|\mathbf{w}\|_{L^{\ell,\infty}})\leq K$ and thus $\displaystyle
\sup_{t>0}t^{\mu}\|\mathcal{Q}(t)\|_{L^{r,\infty}}\leq K$ so that $\mathcal{Q}$ maps $\overline{B}_{K}(0)$ into itself. Arguing as above, we can show without ambiguity that $\mathcal{Q}$ is a contraction map on $\overline{B}_{K}(0)$ for appropriately chosen $K$ (small). Applying the Banach fixed point theorem, we thus obtain the existence of a solution $u$ to (\ref{eq:Duhamel}) in $\mathbf{X}$ which is unique in $\overline{B}_{K}(0)$.
Also, $u(t)\rightarrow u_0$, $t\rightarrow 0^{+}$ in  $\mathcal{S}'(\RN)$. To see this let $\langle\cdot,\cdot\rangle$ denotes the duality bracket between $\mathcal{S}'(\RN)$ and $\mathcal{S}(\RN)$. It suffices to show that
 $ \displaystyle\lim_{t\rightarrow 0^+}|\langle u(t)-u_0,\varphi\rangle|=0$ for all $\varphi\in \mathcal{S}(\RN)$. The justification of this fact seemingly follows the lines of the proof in \cite{G}. The details are therefore omitted.\\
  For the remaining bit of the proof, Part \ref{radial-symmetry} is proved as follows.
The solution constructed above via a fixed point argument can be realized as the limit of the following sequence of approximations
\begin{equation*}
u_1={\mathbf S}_{d,0}u_0,\quad u_{j+1}=u_1+{\mathbf A}_1(t)u_j+W(t,x),\quad j=1,2,3,...
\end{equation*} where
\begin{equation*}{\mathbf A}_1(t)v=\int_{0}^{t}{\mathbf S}_{d,-\alpha}(t-s)|v(s)|^{p}\,ds,\quad W(x,t)=\int_{0}^{t}{\mathbf S}_{d,0}(t-s){\mathbf w}\zeta(s)ds.
\end{equation*}
Assume $u_0, {\mathbf w}$ are radial functions. Since the kernel $E_{d}$ is radial in the $x$-variable, we deduce that $u_1$ is radial in $x\in \mathbb{R}^{N}$ as the convolution of two radial functions. For the same reason, $W(x,t)$ is radial in $x$ for all $0<t<\infty$. Likewise, note that ${\mathbf A}_1(t)v$ is radial in the spatial variable provided $v$ is so. Hence, an induction argument shows that each element of the sequence $(u_j)_{j\geq 1}$ is radial in $x\in \mathbb{R}^{N}$. On the other hand, $u_j$ converges to $u$ in a weak-$\star$ sense in $C_{b}([0,\infty);L^{p_c,\infty}(\mathbb{R}^{N}))$ and up to a subsequence which we still denote by $(u_j)_j$; $u_j\rightarrow u$ almost everywhere as $j\rightarrow \infty$ for all $t\in (0,\infty)$. The conclusion follows from the fact that almost everywhere convergence preserves radial symmetry. The second \ref{monotonicity} and last part \ref{positivity} are established by essentially mimicking the previous argument bearing in mind that radial monotonicity and positivity (under the condition $d\in (0,1]$, $E_{d}$ has a positive kernel) are properties which are preserved under convolution. This achieves the proof of Theorem \ref{theo:global existence}.

\subsection{Proof of Theorem \ref{theo:unconditional-uniqueness}}
The strategy to prove Theorem \ref{theo:unconditional-uniqueness} is exactly the same as that employed in \cite[Theorem 1.1, (ii)]{Tay} and \cite[Theorem 2.4]{G}. Therefore we omit the details and simply refer the interested reader to \cite{BTW, Tay, G}.

\subsection{Proof of Theorem \ref{theo:non-existence}}
Set $p_{F}=p_{F}(m)=\frac{N-2m d+\alpha}{N-2dm-2d}$, $m\in \mathbb{R}$ and let $p\in (1,p_F)$.  Assume $\textbf{w}\in C_0(\mathbb{R}^{N})\cap L^{1}(\mathbb{R}^{N})$ such that $\displaystyle\int_{\mathbb{R}^{N}}\textbf{w}(x)\,dx>0$ and suppose that the function $\zeta$ is given according to \eqref{eq:Z1}. By way of contradiction, assume that Problem (\ref{eq:base}) has a global weak solution in the sense of Definition \ref{defn:weak-solution}. Consider the cut-off functions $\psi_k\in C^{\infty}_{0}([0,\infty))$, $k=1,2$ with $0\leq \psi_k\leq 1$ and
\[\psi_1(s)=\begin{cases}
1 \hspace{0.12cm}\text{if}\hspace{0.12cm} 1/2\leq s\leq 3/4\\
0 \hspace{0.12cm}\text{if}\hspace{0.12cm} s\in [0,1/4]\cup[4/5,\infty)
\end{cases},\quad\psi_2(s)=\begin{cases}
1 \hspace{0.12cm}\text{if}\hspace{0.12cm}  s\in [0,1]\\
0 \hspace{0.12cm}\text{if}\hspace{0.12cm} s\geq 2.
\end{cases}
\]  	
Next, pick $T>0$ large enough and let us introduce the function
\[\psi_{T}(x,t)=\psi_1\bigg(\frac{t}{T}\bigg)^{\frac{p}{p-1}}\psi_2\bigg(\frac{|x|^{2d}}{T}\bigg)^{\frac{2dp}{p-1}}.
\]
Since $u$ is a global weak solution of \eqref{eq:base}, we have that
\begin{align*}
-\int_{\mathbb{R}^{N}}\int_{0}^{T}u\partial_t\psi_T-\int_{\mathbb{R}^{N}}\psi_T(x,0)u_0dx&=\int_{\mathbb{R}^{N}}\int_{0}^{T}u(-\Delta)^{d}\psi_{T}+\int_{\mathbb{R}^{N}}\int_{0}^{T}|x|^{\alpha}|u|^{p}\psi_T\\
&\hspace{3.5cm}\qquad{}+\int_{\RN}\int_{0}^{T}\textbf{w}(x)\zeta(t)\psi_T.
\end{align*}
Remark that
$\displaystyle\int_{\RN}\psi_{T}(x,0)u_0(x)dx=0$ so that the above expression implies
\begin{align}\label{eq:integral-ineq}
\nonumber\int_{\RN}\int_{0}^{T}|x|^{\alpha}|u|^{p}\psi_T
+\int_{\RN}\int_{0}^{T}\textbf{w}(x)\zeta(t)\psi_T&\leq \int_{\RN}\int_{0}^{T}|u(-\Delta)^{d}\psi_{T}|+\\
&\hspace{2.3cm}\qquad{}\int_{\RN\times [0,T]}|u||\partial_t\psi_T|.
\end{align}
At this point, set $I_1=\displaystyle\int_{\mathbb{R}^{N}}\int_{0}^{T}|u(-\Delta)^{d}\psi_{T}|$ and $I_2=\displaystyle \int_{\mathbb{R}^{N}\times [0,T]}|u||\partial_t\psi_T|$. We wish to find suitable bounds for both $I_1$ and $I_2$. Start by observing that via an induction argument, one has
\begin{equation}\label{eq:bound-on-d-Laplace}
|(-\Delta)^{d}\psi_{T}|\leq CT^{-1}\psi_2\bigg(\frac{|x|^{2d}}{T}\bigg)^{\frac{2d}{p-1}}\psi_{1}\left(\frac{t}{T}\right)^{\frac{p}{p-1}}.
\end{equation}
Hence, by utilizing the $\varepsilon$-Young inequality, we arrive at
\begin{equation*}
I_{1}\leq \frac{1}{2}\int_{\RN}\int_{0}^{T}|x|^{\alpha}|u|^{p}\psi_Tdxdt+C\int_{\mathbb{R}^{N}}\int_{0}^{T}|(-\Delta)^{d}\psi_T|^{\frac{p}{p-1}}|x|^{-\frac{\alpha}{p-1}}\psi_T^{-\frac{1}{p-1}}dxdt
\end{equation*} and by invoking \eqref{eq:bound-on-d-Laplace}, the second expression (call it $I_{11}$) in the right hand side of the above estimate can further be estimated as follows
\begin{align*}
I_{11}&\leq
cT^{-\frac{p}{p-1}}\bigg(\int_{0}^{T}\psi_1\bigg(\frac{t}{T}\bigg)^{\frac{p^2}{(p-1)^2}-\frac{p}{(p-1)^{2}}}dt\bigg)\bigg(\int_{\{|y|< 2^{\frac{1}{2d}}\}}T^{-\frac{\alpha}{2d(p-1)}+\frac{N}{2d}}|y|^{-\frac{\alpha}{p-1}}dy\bigg)\\
&\leq cT^{-\frac{p}{p-1}+1}T^{\frac{N}{2d}-\frac{\alpha}{2d(p-1)}}\bigg(\int_{0}^{1}\psi_1(\tau)^{\frac{p}{p-1}}d\tau\bigg)\int_{\{|y|< 2^{\frac{1}{2d}}\}}|y|^{-\frac{\alpha}{p-1}}dy\\
&\leq cT^{-\frac{1}{p-1}-\frac{\alpha}{2d(p-1)}+\frac{N}{2d}}
\end{align*}
where we have made in the first and second lines the change of variable $x=T^{1/2d}y$ and $T\tau=t$, respectively. This bound implies in particular that
\begin{align}\label{eq:bound-on-I1}
I_1&\leq \frac{1}{2}\int_{\mathbb{R}^{n}\times [0,T]}|u|^{p}|x|^{\alpha}\psi_Tdxdt+cT^{-\frac{1}{p-1}-\frac{\alpha}{2d(p-1)}+\frac{N}{2d}}.
\end{align}
Similarly, one has
\begin{align*}
I_2&\leq \frac{1}{2}\int_{\mathbb{R}^{N}\times [0,T]}|u|^{p}|x|^{\alpha}\psi_T+C\int_{\mathbb{R}^{N}\times [0,T]}\bigg|\partial_t\psi_1\bigg(\frac{t}{T}\bigg)^{\frac{p}{p-1}}\bigg|^{\frac{p}{p-1}}\psi_1(t/T)^{-\frac{p}{(p-1)^2}}\psi_{2}^{\frac{2dp}{p-1}}(|x|^{2d}/T)|x|^{-\frac{\alpha}{p-1}}\\
&\leq \frac{1}{2}\int_{\mathbb{R}^{N}\times [0,T]}|u|^{p}|x|^{\alpha}\psi_T+C\bigg(\int_{0}^{T}T^{-\frac{p}{p-1}}|\psi_1'(t/T)|^{\frac{p}{p-1}}dt\bigg) \int_{\mathbb{R}^{N}}|x|^{-\frac{\alpha}{p-1}}\psi_{2}(|x|^{2d}/T)^{\frac{2dp}{p-1}}dx\\
&\leq  \frac{1}{2}\int_{\mathbb{R}^{N}\times [0,T]}|u|^{p}|x|^{\alpha}\psi_T+CT^{\frac{-\alpha}{2d(p-1)}+\frac{N}{2d}-\frac{1}{p-1}}\bigg(\int_{0}^{1}[\psi'_1(s)]^{\frac{p}{p-1}}ds\bigg)\int_{|y| < 2^{\frac{1}{2d}}}|y|^{-\frac{\alpha}{p-1}}\psi_{2}(|y|^{2d})^{\frac{2dp}{p-1}}dy
\end{align*}
so that
\begin{equation}\label{eq:boud-on-I2}
I_2\leq  \frac{1}{2}\int_{\mathbb{R}^{N}\times [0,T]}|u|^{p}|x|^{\alpha}\psi_Tdxdt+CT^{\frac{-\alpha}{2d(p-1)}+\frac{N}{2d}-\frac{1}{p-1}}.
\end{equation}
On the other hand, since $T$ is chosen large we bound the second term in the left-hand side of (\ref{eq:integral-ineq}) from below as follows
\begin{align*}
\int_{\RN\times [0,T]}\textbf{w}(x)\zeta(t)\psi_T(x,t)dxdt&\geq \int_{T/2}^{T}\int_{\mathbb{R}^{N}}\textbf{w}(x)\zeta(t)\psi_T(x,t)dxdt\\
&\geq \bigg(\int_{T/2}^{T}t^{m}\psi_1\bigg(\frac{t}{T}\bigg)^{\frac{p}{p-1}}dt\bigg)\bigg(\int_{\mathbb{R}^{N}}\textbf{w}(x)\psi_2\bigg(\frac{|x|^{2d}}{T}\bigg)^{\frac{2d p}{p-1}}dx\bigg)\\
&\geq T^{m+1}\bigg(\int_{1/2}^{1}\psi_1(\tau)^{\frac{p}{p-1}}d\tau\bigg)\int_{\mathbb{R}^{N}}\textbf{w}(x)\psi_2\bigg(\frac{|x|^{2d}}{T}\bigg)^{\frac{2d p}{p-1}}dx.
\end{align*}
From (\ref{eq:integral-ineq}) and by combining the above bound with (\ref{eq:bound-on-I1}) and (\ref{eq:boud-on-I2}), we obtain that
\[\int_{\RN}\textbf{w}(x)\psi_2\bigg(\frac{|x|^{2d}}{T}\bigg)^{\frac{2dp}{p-1}}dx\leq CT^{-\frac{\alpha}{2d(p-1)}+\frac{N}{2d}-\frac{p}{p-1}-m}.
\]
We deduce from the latter that
\[\int_{\RN}\textbf{w}(x)dx\leq 0
\]
which yields a contradiction in view of the assumption imposed on $\textbf{w}$ along with the condition $p<p_{F}$. In fact, the function $\psi_2\bigg(\frac{|x|^{2d}}{T}\bigg)^{\frac{2dp}{p-1}}$ converges in a pointwise sense to $1$ as $T$ approaches infinity so that by applying the Dominated Convergence Theorem, there holds
\[\lim_{T\rightarrow \infty} \int_{\RN}\textbf{w}(x)\psi_2\bigg(\frac{|x|^{2d}}{T}\bigg)^{\frac{2d p}{p-1}}dx=\int_{\RN}\textbf{w}(x)dx>0.
\]
As a consequence, for sufficiently large $T$, $\displaystyle\int_{\RN}\textbf{w}(x)\psi_2\bigg(\frac{|x|^{2d}}{T}\bigg)^{\frac{2dp}{p-1}}dx\geq \gamma\int_{\RN}\textbf{w}(x)dx$ for some $\gamma<1$. This achieves the proof of the first part of Theorem \ref{theo:non-existence}. When $m>0$, one can easily adapt the preceding approach with a slightly modified test function. More precisely, for $R>0$, by replacing $\psi_T$ by the function $\psi_{T,R}(x,t)=\psi_1\big(\frac{t}{T}\big)^{\frac{p}{p-1}}\psi_2\big(|x|^{2d}R^{-2d}\big)^{\frac{2dp}{p-1}}$, we reach the same conclusion $\displaystyle\int_{\RN}\textbf{w}(x)dx\leq 0$ which leads to a contradiction. Theorem \ref{theo:non-existence} is now completely proved.
\section*{ \bf Acknowledgements}
\noindent The author would like to express his deep thanks to {\em Gael Diebou Yomgne} for interesting discussions. He also thanks {\em Slim Tayachi} for helpful remarks and suggestions.

\noindent$\rule[0.05cm]{14cm}{0.05cm}$

	\medskip

\noindent{\bf\large Funding.}
Funding information is not applicable / No funding was received. 

\noindent{\bf\large Declarations.}
On behalf of all authors, the corresponding author states that there is no conflict of interest. 
No data-sets were generated or analyzed during the current study.

\medskip
\noindent$\rule[0.05cm]{14cm}{0.05cm}$



\begin{thebibliography}{10000}

\bibitem{BTW}{B. Ben Slimene, S. Tayachi and F. B. Weissler}, {\em Well-posedness, global existence and large time behavior for Hardy-H\'enon parabolic equations}, Nonlinear Anal., {\bf 152} (2017), 116--148.

\bibitem{BLZ} C. Bandle, H. A. Levine and Qi S. Zhang, {\em Critical Exponents of Fujita Type for Inhomogeneous Parabolic Equations and Systems}, Journ. of Math. Anal. and App., {\bf 251} (2000), 624--648.

\bibitem{Caz-Dick-Weiss} T. Cazenave, F. Dickstein and Fred B. Weissler, {\em An equation whose Fujita critical exponent is not given by scaling}, Nonlinear  Analysis, {\bf 68} (2008), 862--874.

\bibitem{fujita}{H. Fujita}, {\em On the blowing up of solutions of the Cauchy problem for $u_t=\Delta u+u^{1+\alpha}$}, J. Fac. Sci. Univ. Tokyo Sec. IA Math., {\bf 13} (1966), 109--124.

\bibitem{GP} {V. A. Galaktionov  and   S. I. Pohozaev }, {\em Existence and blow-up for higher-order semilinear parabolic equations: majorizing order-preserving operators}, {Indiana Univ. Math. J.}, {\bf 51} (2002), 1321--1338.

\bibitem{GV} {V. A. Galaktionov and J. L. V\'azquez}, {\em The problem of blow-up in nonlinear parabolic equations}, Discrete Contin. Dyn. Syst., {\bf 8} (2002), 399--433.

\bibitem{Hayak} {K. Hayakawa}, {\em On nonexistence of global solutions of some semilinear parabolic differential equations}, Proc. Japan Acad., {\bf 49} (1973), 503--505.

 \bibitem{Hu} {B. Hu},  {\em Blow Up Theories for Semilinear Parabolic Equations}, Springer, Berlin (2011).


\bibitem{Jacob} N. Jacob, {\em Pseudo-differential operators and Markov processes. Vol. I. Fourier analysis and semigroups}, Imperial College Press, London, (2001).

\bibitem{JKS} M. Jleli, T. Kawakami and B. Samet, {\em Critical behavior for a semilinear parabolic equation with forcing term depending of time and space},  J. Math. Anal. Appl., {\bf 486} (2020), 123931.

\bibitem{LSU} O.A. Ladyzenskaja, V.A. Solonnikov and N.N. Ural\'ceva, {\em Linear and quasilinear equations of parabolic type}, Amer. Math. Soc., Transl. Math. Monographs, Providence, R.I.(1968).

\bibitem{Lev} H. A. Levine, {\em The role of critical exponent in blow-up theorems}, SIAM. Rev., {\bf 32} (1990), 262--288.

\bibitem{M} M. Kwasnicki, {\em Ten equivalent definitions of the fractional Laplace operator}, Fract. Calc. Appl. Anal., {\bf 20}(1) (2017), 7--51.


\bibitem{PGL} {P. G. Lemari{\'e}-Rieusset}, {\em Recent developments in the {Navier}-{Stokes} problem}, {Chapman Hall/CRC Res. Notes Math.}, Vol. {\bf 431}, {2002}.

\bibitem{Ma} M. Majdoub, {{\em Well-posedness and blow-up for an inhomogeneous semilinear parabolic equation}}, {{Differ. Equ. Appl.}}, {\bf 13} (2021), {85--100}.


 \bibitem{ADE} M. Majdoub, S. Otsmane and S. Tayachi, {\em Local Well-posedness and Global Existence for the Biharmonic Heat Equation with Exponential Nonlinearity}, Advances in Differential Equations, {\bf 23} (2018), 489--522.

 \bibitem{MT1} {A. V. Martynenko and A. F. Tedeev}, {\em Cauchy problem for a quasilinear parabolic equation with a source term and an inhomogeneous density}, Comput. Math. Math. Phys., {\bf 47} (2007), 238--248.

\bibitem{MT2} {A. V. Martynenko and A. F. Tedeev}, {\em On the behavior of solutions to the Cauchy problem for a degenerate parabolic equation with inhomogeneous density and a source}, Comput. Math. Math. Phys., {\bf 48} (2008), 1145--1160.


 \bibitem{MP} E. Mitidieri and S. I. Pohozaev, {\em A priori estimates and blow-up of solutions of nonlinear partial differential equations and inequalities}, Proc. Steklov Inst. Math., {\bf 234} (2001), 3--383.

\bibitem{Qi} Yuan-wei Qi, {\em The critical exponents of parabolic equations and blow-up in $\R^n$}, Proceedings of the Royal Society of Edinburgh, {\bf 128A} (1998), 123--136.

\bibitem{QS} {P. Quittner and P. Souplet}, {\em Superlinear parabolic problems}, {Birkh\"auser Verlag, Basel} (2007), {xii+584}.

\bibitem{RT} {B. Ruf and E. Terraneo}, {\em The {Cauchy} problem for a semilinear heat equation with singular initial data}, {Evolution equations, semigroups and functional analysis. In memory of Brunello Terreni. Containing papers of the conference, Milano, Italy, September 27-28, 2000}, {295--309}, {2002}.


\bibitem{Tay} S. Tayachi, {\em Uniqueness and non-uniqueness of solutions for critical Hardy-H\'enon parabolic equations}, J. Math. Anal. Appl., {\bf 488} (2020), 123976.

\bibitem{Terr1} {E. Terraneo}, {\em On the non-uniqueness of weak solutions of the nonlinear heat equation with nonlinearity {{\(u^3\)}}}, {C. R. Acad. Sci., Paris, S{\'e}r. I, Math.}, {\bf 328} (1999), {759--762}.,
 Year = {1999},

\bibitem{Terr2}{E. Terraneo}, {\em Non-uniqueness for a critical nonlinear heat equation}, {Commun. Partial Differ. Equations}, {\bf 27} (2002), {185--218}.


 \bibitem{Weissler}{ F. B. Weissler},
{\em Semilinear evolution equations in Banach spaces},  J. Funct. Anal.,  {\bf 32} (1979), no. 3, 277--296.

\bibitem{Ya}{ M. Yamazaki},
{\em The Navier-Stokes equations in the weak-$L^{n}$ space with time-dependent external force},  Math. Ann., {\bf 317} (2000), 635--675.

\bibitem{G}{G. D. Yomgne}, {\em On the generalized parabolic Hardy-Henon equation: Existence, blow-up, self-similarity and large-time asymptotic behavior}, Diff. Int. Equ., {\bf 35} (2022), 57--88.

\bibitem{Zh1} Q. S. Zhang, {\em A new critical phenomenon for semilinear parabolic problem}, J. Math. Anal. Appl., {\bf 219} (1998), 123--139.

\bibitem{Zh2} Q. S. Zhang, {\em Blow up and global existence of solutions to an inhomogeneous parabolic system}, J. Differential Equations, {\bf 147} (1998), 155--183.

\bibitem{Wang} X. Wang, {\em On the Cauchy problem for reaction-diffusion equations}, Trans. Amer. Math. Soc., 337 (1993), 549--590.



\end{thebibliography}
\end{document}